\documentclass{amsart}

\usepackage{amssymb}
\usepackage[all]{xy}
\newtheorem{thm}{Theorem}[section]

\newtheorem{conj}[thm]{Conjecture}
\newtheorem{main-conj}[thm]{Main Conjecture}
\newtheorem{main-ques}[thm]{Main Question}

\theoremstyle{definition}
\newtheorem{defn}[thm]{Definition}

\newtheorem{say}[thm]{}
\newtheorem{exmp}[thm]{Example}

\newtheorem*{ack}{Acknowledgments}      % \renewcommand{\theack}{} 

\newtheorem{defn-thm}[thm]{Definition--Theorem}  %!!!!!!!!!!!!!!!!!!!!!!!!
\newtheorem{defn-lem}[thm]{Definition--Lemma}  %!!!!!!!!!!!!!!!!!!!!!!!!
\newtheorem{brot}[thm]{Basic rule of thumb}

\theoremstyle{remark}

\let \cedilla =\c

\newcommand{\bb}[0]{{\mathbb B}}
\newcommand{\zz}[0]{{\mathbb Z}}

\newcommand{\rr}[0]{{\mathbb R}}

\newcommand{\cc}[0]{{\mathbb C}}

\newcommand{\pp}[0]{{\mathbb P}}
\newcommand{\qq}[0]{{\mathbb Q}}
\newcommand{\map}[0]{\dasharrow}
\newcommand{\qtq}[1]{\quad\mbox{#1}\quad}

\newcommand{\sing}[0]{\operatorname{Sing}}

\def\into{\DOTSB\lhook\joinrel\to}

\def\loccoh#1.#2.#3.#4.{H^{#1}_{#2}(#3,#4)}

\DeclareMathAlphabet{\mathchanc}{OT1}{pzc}%
                                {m}{it}

\newcommand{\GL}{\mathrm{GL}}

\title{The structure of algebraic varieties}
\author{J\'anos Koll\'ar}

\begin{document}
\bibliographystyle{amsplain}   %{amsplain} %{amsalpha}

\maketitle

\begin{abstract}
The aim of this address is to give an overview of the main
questions and results of the structure theory of higher dimensional 
algebraic varieties.
\end{abstract}

%\keywords{algebraic variety, Mori program, moduli questions}

%\subjclass{14-02, 14E30, 14B05,  14D20   }

\section{Early history: Euler, Abel, Jacobi, Riemann}

Our story, like many others in mathematics, can be traced back at least
 to Euler who studied elliptic integrals of the form
$$
\int \frac{dx}{\sqrt{x^3+ax^2+bx+c}}.
$$
The study of integrals of algebraic functions was further developed
by Abel and Jacobi. From our point of view the next major step
was taken by Riemann. Instead of dealing with a multi-valued function
like $\sqrt{x^3+ax^2+bx+c} $, Riemann looks at the 
complex algebraic curve
$$
C:=\bigl\{(x,y): y^2=x^3+ax^2+bx+c\bigr\}\subset \cc^2.
$$
Then the above integral becomes
$$
\int_{\Gamma} \frac{dx}{y}
$$
for some path $\Gamma$  on the  
algebraic curve $C$.
More generally, a polynomial $g(x,y)$ implicitly defines
 $y:=y(x)$ as a multi-valued function of $x$ and
for any  meromorphic function $h(u,v)$, the multi-valued integral
$$
\int h\bigl(x, y(x)\bigr)\ dx
$$
becomes a single valued integral 
$$
\int_{\Gamma} h\bigl(x, y\bigr)\ dx
$$
for some path $\Gamma$ on the  
algebraic curve
$C(g):=\bigl(g(x,y)=0\bigr)\subset \cc^2$.
Substitutions that transform one integral associated
to a polynomial $g_1$  into another integral associated
to a  $g_2$
can be now seen as algebraic maps between the 
curves $C(g_1) $ and  $C(g_2) $.

Riemann also went  further. As a simple example, consider
the curve $C$ defined by  $(y^2=x^3+x^2)$ and notice that
$(t^3-t)^2\equiv(t^2-1)^3+(t^2-1)^2$.
Thus the substitution
$x=t^2-1,\ y=t^3-t$ (with inverse $t=y/x$) allows us to transform any integral
$$
\int h\bigl(x, \sqrt{x^3+x^2}\bigr)\ dx
\qtq{into}
\int h\bigl(t^2-1, t^3-t\bigr)\cdot 2tdt.
$$
To put it somewhat differently, the map
$$
t\mapsto  \bigl(x=t^2-1, y=t^3-t\bigr)
\qtq{and its inverse} 
(x,y)\mapsto  t=y/x
$$
establish an isomorphism
$$
\left\{
\begin{array}{c}
\mbox{meromorphic functions}\\
\mbox{on the curve $(y^2=x^3+x^2)$}
\end{array}
\right\}
\leftrightarrow
\left\{
\begin{array}{c}
\mbox{meromorphic functions}\\
\mbox{on the complex plane $\cc$}
\end{array}
\right\}.
$$
It is best to work with meromorphic functions on $\cc$ that are
also meromorphic at infinity; these live naturally on the
Riemann sphere  $\cc\pp^1$. We can now state Riemann's
fundamental theorem as follows.

\begin{thm}[Riemann, 1851] \label{riemann.thm}
For every algebraic curve
$C\subset \cc^2$ there is a unique, compact Riemann surface $S$
and a meromorphic map  $\phi:S\map C$ with meromorphic inverse
$\phi^{-1}:C\map S$ such that
$$
f_C\mapsto  f_S:=f_C\circ \phi\qtq{and} f_S\mapsto  f_C:=f_S\circ \phi^{-1}
$$
establish an isomorphism between the meromorphic function theory on
$C$ and the meromorphic function theory on
$S$.
\end{thm}

\section{Main questions, informally}

We can now give an initial formulation of the two main problems
that we  consider; the precise versions are stated in
Sections \ref{blocks.sec} and \ref{moduli.sec}.
The first is a direct higher-dimensional analog of the 
results of Riemann. (See Section \ref{sec.defn} for basic definitions.)

\begin{main-ques} \label{prob1}
Given an algebraic variety $X$, is there another algebraic
 variety $X^{\rm m}$ such that
\begin{enumerate}
\item the meromorphic  function theories of $X$ and of $X^{\rm m}$
are isomorphic and
\item the geometry of $X^{\rm m}$ is the ``simplest'' possible?
\end{enumerate}
\end{main-ques}

Riemann's theorem says that, in dimension 1, ``simplest''
should mean smooth and compact, but  in higher dimensions
smoothness is not the right notion. One of the hardest aspects
of the theory was to understand what the correct concept
of ``simplest'' should be.

So far we have dealt with individual algebraic varieties.
A salient feature of algebraic geometry is that by
continuously varying the coefficients of the defining polynomials
we get continuously varying families of  algebraic varieties.
We can thus study how to transform a family $\{X_t:t\in T\}$
of varieties into its ``simplest'' form.
A tempting idea is to take the ``simplest'' forms
$\{X_t^{\rm m}:t\in T\}$ obtained previously.
Unfortunately, this fails already in dimension 1.
Starting with a family of curves $\{C_t:t\in T\}$,
the corresponding  Riemann surfaces $\{S_t:t\in T\}$
form a  continuously varying family over a dense
open subset $T^0\subset T$ but not everywhere.

For curves the correct answer was found by Deligne and Mumford in 1969.
We  use the guidance provided by this 1-dimensional case
and the answer to the first Main Question to answer the
second.

\begin{main-ques}  \label{prob2}{\ }
 What are the ``simplest'' families of algebraic varieties?
How can one transform an arbitrary family into one of  the  ``simplest'' 
families?
% Is there a good way to describe all  the  ``simplest'' families?
\end{main-ques}

\section{What are algebraic varieties?}\label{sec.defn}

Here we quickly recall the basic concepts and definitions that we use.
For general introductory texts, see \cite{shaf, MR0453732, Harris95}.

An {\it affine algebraic set} in   $\cc^N$ is
the common zero-set of  some polynomials
$$
\begin{array}{rcl}
X^{\rm aff}&=&X^{\rm aff}(f_1,\dots, f_r)\\
&=&\bigl\{(x_1,\dots, x_N): 
f_1(x_1,\dots, x_N)=\cdots=f_r(x_1,\dots, x_N)=0\bigr\}\subset \cc^N.
\end{array}
$$
It is especially easy to visualize {\it hypersurfaces}  $X(f)\subset \cc^N$
 defined by 1 equation.
Usually we  count complex dimensions,
thus $\dim \cc^N=N$ and  $\dim X$ is one half of the
usual topological dimension of $X$. In low dimensions we talk about
{\it curves, surfaces, 3-folds.} 
Thus, somewhat confusingly, an algebraic curve is  a
(possibly singular) Riemann surface.

An affine algebraic set $X$ is called {\it irreducible}
if it can not be written as a union of two algebraic sets in a nontrivial way.
 Such sets are
called {\it affine algebraic varieties.} Every algebraic set $X$ is a
finite union of  algebraic varieties  $X=\cup_iX_i$ such that
$X_i\nsubseteq X_j$ for $i\neq j$. Such a decomposition is unique,
up to permuting the indices. 
Thus from now on we are interested mainly
in algebraic varieties.

For example, the irreducible components of a hypersurface  $X(f)$
correspond to the  irreducible factors of $f$, thus $X(f)$
 is irreducible iff $f$ is a power of an
irreducible polynomial.

An affine algebraic set $X^{\rm aff} $ is compact iff it is 0-dimensional,
thus it is almost always better to work with the closure
of $X^{\rm aff}$ in the complex projective space
$$
X:=X^{\rm proj}\subset \cc\pp^N.
$$
Thus we get  {\it projective algebraic sets} and 
{\it projective varieties.} Finally, a {\it quasi-projective variety}
is an open subset $U$ of a projective variety $X$ whose complement
$X\setminus U$ is  a projective algebraic set. Note that $U$ is a 
``very large'' subset of $X$, in particular $U$ is dense in $X$.
This is a key feature of algebraic geometry: all open subsets are
 ``very large.''

On a complex projective space $\cc\pp^N$ the
homogeneous coordinates $(x_0{:}\cdots{:}x_N)$ 
are defined only up to multiplication by a scalar.
Thus one can not evaluate a  polynomial 
$p(x_0,\dots, x_N)\in \cc[x_0,\dots, x_N]$,
at a point of  $\cc\pp^N$. However, if $p$ is homogeneous of degree $d$
then
$$
p(\lambda x_0,\dots, \lambda x_N)=\lambda^dp(x_0,\dots, x_N).
$$
Thus the zero set of $p$ is well-defined and a
quotient of two  homogeneous polynomials of the same degree
$$
f(x_0,\dots, x_N)=\frac{p_1(x_0,\dots, x_N)}{p_2(x_0,\dots, x_N)}
$$
is also well-defined (except where $p_2$ vanishes). 
These are the {\it rational functions}
 on $\cc\pp^N$. By restriction, we get rational functions
on any projective variety $X\subset \cc\pp^N$.

At first sight these seem downright antiquated definitions; a modern
theory ought to  be local. That is, one should 
consider varieties that are locally defined by analytic functions
and work with  meromorphic  functions on them. However, we know that
every  meromorphic  function on $\cc\pp^1$ is rational and 
the same holds in all dimensions.

\begin{thm}[Chow, 1949; Serre, 1956]
 Let $M\subset \cc\pp^N$ be a closed subset that
can be {\em locally}  given as the common zero set of analytic functions. Then
\begin{enumerate}
\item $M$ is algebraic, that is, it can be {\em globally} given as 
the common zero set of   homogeneous polynomials and
\item every meromorphic  function $f$ on $M$ is algebraic, 
that is, $f$ can be {\em globally} given as 
the quotient of two  homogeneous polynomials of the same degree.
\end{enumerate}
\end{thm}

Now we come to a key feature of algebraic geometry.
There are two competing notions of ``map'' and
two competing notions of ``isomorphism.''

\begin{defn}[Map and morphism] Let $X\subset \cc\pp^N$ be an algebraic variety
and $f_0,\dots, f_M$ nonzero rational functions on $X$.
They define a {\it  map}  (or rational map)
$$
{\mathbf f}:X\map \cc\pp^M
\qtq{given by} p\mapsto \bigl(f_0(p){:}\cdots {:} f_M(p)\bigr)\in \cc\pp^M.
$$
To start with, ${\mathbf f }$ is only defined 
at a point $p$ if
 none of the $f_i$ has a pole at $p$ and
 not all of the $f_i$ vanish at $p$.  However, 
since the projective coordinates are defined only up to a scalar multiple,
$(gf_0,\dots, gf_M)$ define the same map for any rational function $g$,
thus it can happen that  ${\mathbf f }$ is everywhere defined.
In this case it is called a {\it morphism.}
A map is denoted by $\map$ and a morphism by $\to$.

For example, projecting $\cc\pp^2$ from the origin  $(0{:}0{:}1)$ to the
line at infinity is given by  
$$
\pi:(x{:}y{:}z)\mapsto  \bigl(\tfrac{x}{z} : \tfrac{y}{z}\bigr)=
\bigl(\tfrac{x}{y} : 1\bigr)=\bigl(1 : \tfrac{y}{x}\bigr).
$$
Thus $\pi$ is  defined everywhere except at  $(0{:}0{:}1)$.
\end{defn}

\begin{defn}[Isomorphism]
Two quasi-projective varieties  $X\subset \cc\pp^N$ and  $Y\subset \cc\pp^M$
are {\it isomorphic} if there are morphisms
$$
f:X\to Y \qtq{and}  g:Y\to X
$$
that are inverses of each other. Isomorphism is denoted by $X\cong Y$.
\end{defn}

We will think of isomorphic varieties as being essentially the same.
Using maps instead of morphisms in the above definition yields
the notion of birational equivalence. This notion is
unique to algebraic geometry; it has no known analog in
topology or differential geometry.

\begin{defn}[Birational equivalence]\label{bir.eq.defn}
Two  quasi-projective varieties  $X\subset \cc\pp^N$ and  $Y\subset \cc\pp^M$
are {\it birational} 
(in old terminology, birationally isomorphic)
if there are rational maps
$$
f:X\map Y \qtq{and}  g:Y\map X
$$
such that the following equivalent conditions hold.
\begin{enumerate}
\item $\phi_Y\mapsto  \phi_X:=\phi_Y\circ f$ and 
$\phi_X\mapsto  \phi_Y:=\phi_X\circ g$
establish an isomorphism between the meromorphic (=rational) function theory on
$X$ and the meromorphic (=rational) function theory on $Y$.
\item There are algebraic subsets $Z\subsetneq X$ 
and $W\subsetneq Y$ such that
$(X\setminus Z)\cong (Y\setminus W)$. 
%$f|_{X\setminus Z}$ and $g|_{Y\setminus W}$ are everywhere defined
%and  inverses of each other.
\end{enumerate}

As an example, consider the affine surface $S:=(x^2+y^2=z^3)\subset \cc^3$.
It is birational to $\cc^2_{uv}$ as shown by the rational maps
$$
f:(x,y,z)\map \bigl(\tfrac{x}{z}, \tfrac{y}{z}\bigr)
\qtq{and}
g:(u,v)\to \bigl(u(u^2+v^2), v(u^2+v^2), u^2+v^2\bigr).
$$
Here $f$ is not defined if $z=0$ while $g$ is everywhere defined
but it maps the pair of lines  $(u=\pm iv)$ to the 
origin  $(0,0,0)$. Thus
$$
S\setminus(z=0)\cong \cc^2\setminus(u^2+v^2=0)
\qtq{but}  S\not\cong \cc^2.
$$
\end{defn}

\begin{brot} Let $X, Y$ be  algebraic varieties that are
birational to each other.
Many questions of algebraic geometry about $X$
can be answered by 
\begin{itemize}
\item first studying the same question on $Y$ and then 
\item studying a similar question involving the
lower dimensional algebraic sets $Z$ and $W$ as in (\ref{bir.eq.defn}.2).
\end{itemize}
The aim of the Minimal Model Program is to exploit this
in two steps.
\begin{itemize}
\item Given a question  and a variety $X$,
find a  variety $Y$ that is
birational to $X$ such that the 
geometry of $Y$ is ``best adapted'' to studying
the particular question. This is a variant of the  first Main Question.
% \ref{prob1}.
\item Set up the appropriate dimension induction to deal with
the exceptional sets $Z\subset X$ and $W\subset Y$.
\end{itemize}
\end{brot}

{\it Important aside.}
  More generally, if we decompose an algebraic variety into
disjoint locally closed pieces, then the collection of the
pieces carries a lot of information about the variety.
I would like to stress that this is  a rather 
noteworthy fact about algebraic geometry. For instance, if we decompose a
simplicial complex into its simplices, then usually the only information
we retain is the dimension and the Euler characteristic.
By contrast, all the homology groups of a smooth, projective algebraic 
 variety can be recovered from
the pieces. This is a key consequence of Hodge Theory, as
formulated by Deligne, and is a starting point of 
Grothedieck's theory of motives.

\section{Classical results}

After the study of algebraic curves, two main avenues of
investigations were pursued. One direction focused on the local study
of varieties with a main aim of resolving them completely.
The other direction aimed to understand the global structure
of algebraic surfaces. These are both still very active research areas.
We recall a few of the main results that are relevant for the
general theory. For  detailed treatments and for
references see \cite{bpv, k-res}.

\subsection*{Resolution of singularities} {\ }

Riemann's theorem says that
every singular algebraic curve $C$ is birational to a smooth, compact curve
(or Riemann surface). The first steps toward answering the Main Questions
in higher dimensions focused on this problem: {\it Is every algebraic variety
 birational to a smooth, projective variety?}

\begin{defn} A variety $X\subset \cc^N$ is smooth and has dimension $d$
at a point $p\in X$ iff the following equivalent conditions hold.
\begin{enumerate}
\item $X\subset \cc^N\cong\rr^{2N}$ is a $C^{\infty}$-submanifold
of (real) dimension $2d$ near $p$.
\item One can choose coordinates $z_1,\dots, z_N$ and
equations $f_1,\dots, f_{N-d}$ of $X$ such that 
 $(f_1=\cdots = f_{N-d}=0)$ coincides with $X$ near $p$ and
 the Jacobian matrix  
$\bigl(\partial f_i/\partial z_j: 1\leq i,j\leq N-d\bigr)$ is invertible at $p$.
\item There are holomorphic functions $\phi_i=\phi_i(w_1,\dots, w_d)$
defined near the origin  and constants $c_i$ such that
$$
(w_1,\dots, w_d)\mapsto 
\bigl(\phi_1({\bf w}), \dots, \phi_{N-d}({\bf w}), w_1+c_1, \dots, w_d+c_d\bigr)
$$
maps a small {\it ball} 
${\mathbf 0}\in \bb^d(\epsilon)\subset \cc^d$ onto a neighborhood of $p\in X$. 
\end{enumerate}
In the latter case we view $(w_1,\dots, w_d) $ as
{\it local analytic coordinates} on $X$ near $p$. 
(It is an ever present technical problem that there is
no good notion of local algebraic coordinates. Open algebraic neighborhoods
are too large to admit a single-valued coordinate system.)
\end{defn}

On an algebraic variety $X$ the set of singular points turns out 
to be an algebraic subset, denoted by  $\sing X\subset X$.
For every variety $X$, a generalization of Riemann's method (\ref{riemann.thm})
 produces a new variety  $X^{\rm n}\to X$,
called the {\it normalization} of $X$,  such that
$\sing \bigl(X^{\rm n}\bigr)$ has codimension $\geq 2$ in $X^{\rm n}$.
Thus, in higher dimensions, one usually works with
{\it normal varieties} whose singular set has codimension $\geq 2$.

To make the singular set even smaller, or to get rid of it
completely, turned out to be very difficult. The final result
was established by Hironaka in 1964.

\begin{thm}[Resolution of singularities]\label{res.sing.thm}
 For every algebraic variety $X$, there are  (very many)
  smooth, projective varieties $X^{\rm sm}$
birational to $X$. 

If $X$ is projective, one can arrange to have
 a morphism  $f:X^{\rm sm}\to X$ that is an isomorphism
over $X\setminus \sing X$. 
\end{thm}

\subsection*{Algebraic surfaces} {\ }

By resolution of singularities,  any projective surface $S$
is birational to a smooth projective surface $S^{\rm sm}$, but,
in contrast with the theory of curves, there are many such
smooth projective surfaces $S^{\rm sm}$. 
We can thus reformulate the first Main Question: % \ref{prob1}:
 Is there a ``simplest'' one 
among all smooth projective surfaces birational to $S$?

To answer this question, first we study how to make a
smooth projective surface more ``complicated.''

\begin{defn}[Blowing-up] Let $S$ be a smooth  algebraic surface
and $p\in S$ a  point. {\it Blowing-up} is an operation
that creates a new smooth surface $B_pS$ by removing $p$ and replacing it
with a $\cc\pp^1$ corresponding to all the tangent directions 
of $S$ at $p$. Collapsing the new $\cc\pp^1$ to a point gives
a morphism $\pi:B_pS\to S$. 

In local coordinates, it can be described as follows.
Start with the unit ball $\bb^2_{xy}\subset \cc^2_{xy}$
and  $\cc\pp^1_{st}$ where the subscripts name the  coordinates.  Set
$$
B_{\mathbf 0}\bb^2_{xy}:=(xt-ys=0)\subset \bb^2_{xy}\times \cc\pp^1_{st}.
$$
Let  $\pi:B_{\mathbf 0}\bb^2_{xy}\to  \bb^2_{xy}$
denote the coordinate projection.
 If $(x,y)\neq (0,0)$ then $\pi^{-1}(x,y)$
is the single point  $(x,y)\times(x{:}y)$. However, if
$x=y=0$ then  $(s{:}t)$ can be arbitrary, thus 
$\pi^{-1}(0,0)\cong \cc\pp^1_{st}$. 

Note that $s/t=x/y$ is the natural coordinate on $\cc\pp^1_{st}$,
thus blowing up is akin to switching to polar coordinates
since the polar angle $\theta$ equals $\tan^{-1}(x/y)$.
\end{defn}

One can blow up any number of points of $S$ and then repeat
by blowing up some of the new points of $B_pS$. Thus
blowing up is a cheap way to get infinitely many new smooth surfaces out of one.

\begin{defn} A smooth projective surface is called {\it minimal}
if it can not be obtained from another smooth, projective surface
by blowing up.
\end{defn}

This notion allows us to get a very good analog of Riemann's theorem
\ref{riemann.thm}.

\begin{thm}[Enriques, 1914 ;  Kodaira, 1966]\label{min.mod.surf.thm}
For every   projective, algebraic surface $S$, 
exactly one of the following holds.
\begin{enumerate}
\item (Minimal model) There is a unique, minimal surface $S^{\rm m}$
birational to $S$.
\item  $S$ is birational to $C\times \cc\pp^1$
for a unique,   smooth, projective curve $C$.
\end{enumerate}
\end{thm}

\begin{say}[Du~Val singularities]\label{duval.say}
It was also gradually understood that instead of working with the
minimal model $S^{\rm m}$, it is sometimes better to use a
slightly singular {\it canonical model} $S^{\rm can}$.
The resulting singularities were first classified by 
Du~Val in 1934; the list is quite short, ranging from the simplest
$(x^2+y^2+z^2=0)$ to the most complicated $(x^2+y^3+z^5=0)$.
They are also called  {\it rational double points.}

Their importance was not generally
recognized until the 1960's when they were rediscovered from many
different points of view;
see \cite{durfee} for a
survey.
\end{say}

\section{The  first Chern class and the Ricci curvature}

The  first Chern class,
which is closely related to the Ricci curvature,
carries much of the important information about the
structure of a variety. We follow the differential geometry sign conventions;
algebraic geometers usually work with the {\it canonical class,}
which is (a slight refinement of)  the negative of the first Chern class.

\begin{say}[Complex volume forms] A measure on $\rr^n$ can be 
identified with  an $n$-form
$$
s(x_1,\dots, x_n)\cdot dx_1\wedge\dots\wedge dx_n.
$$
Thus a measure on a real manifold $M$ is an $n$-form 
that in local coordinates can be written as above.

Similarly, on a smooth variety $X$ of dimension $n$
a {\it complex volume form}
is an $n$-form   $\omega$ 
that in local holomorphic coordinates  can be written as 
$$
h(z_1,\dots, z_n)\cdot dz_1\wedge\dots\wedge dz_n.
$$
Thus a complex volume form  $\omega$  gives a real volume form
$\left(\frac{\sqrt{-1}}{2}\right)^n\omega\wedge\bar\omega$
where the constant comes from the formula
$$
dz\wedge d\bar z=(dx+\sqrt{-1}dy)\wedge(dx-\sqrt{-1}dy)=
-2\sqrt{-1}\ dx\wedge dy.
$$
(There is usually an additional $\pm$, depending on one's
orientation conventions.)
%and we use the plus sign iff $n\equiv 0,1\mod 4$.
\end{say}

From the point of view of differential geometry,
one would like to use  $C^{\infty}$ complex volume forms,
that is,  the  $h(z_1,\dots, z_n) $ should be  nowhere zero $C^{\infty}$-functions.
Algebraic geometry, however, prefers meromorphic  volume forms
 where the $h(z_1,\dots, z_n) $ are meromorphic functions.
(See (\ref{jac.sing}.1) for some explicit examples.)
Thus the ideal situation is when a complex volume form
is  given by  nowhere zero holomorphic functions
$h(z_1,\dots, z_n) $. This is possible only for 
Calabi--Yau 
varieties; they form a very special but important subclass
(\ref{blocks.ssc.2}). 

Thus in general we try to understand how to connect
$C^{\infty}$ and  meromorphic volume forms. 

On the differential geometry side the key notion is the
curvature which defines the Chern form.

\begin{defn} [Chern form and Chern class]
Let $\omega$ be a  $C^{\infty}$ complex volume form.
%Its {\it curvature form} is the 2-form given in local coordinates as
%$$
%\Theta(\omega)=-2\partial\bar\partial \log |h(z_1,\dots, z_n)|=
%-2\sum_{ij}\frac{\partial^2 \log |h({\mathbf z})|}
%{\partial z_i\partial \bar z_j}dz_i\wedge d\bar z_j.
%$$
The {\it first Chern form}  or {\it Ricci curvature form} of $(X,\omega)$
 is the  2-form
$$
\tilde c_1(X,\omega):= \frac{\sqrt{-1}}{\pi}
\partial\bar\partial \log |h(z_1,\dots, z_n)|=
\frac{\sqrt{-1}}{\pi}\sum_{ij}\frac{\partial^2 \log |h({\mathbf z})|}
{\partial z_i\partial \bar z_j}dz_i\wedge d\bar z_j.
$$
As a 2--form, this depends on the choice of the  volume form $\omega$, but
%essentially by Theorem \ref{gauss-b.thm}, 
it gives a
well defined  De~Rham cohomology class $c^{\rr}_1(X)\in H_{DR}^2(X, \rr)$
which actually lifts to an  integral cohomology class 
$$
c_1(X)\in H^2(X, \zz),
$$
called the {\it first Chern class} of $X$.
\end{defn}
 
\begin{defn}[Algebraic degree] 
  Let $X$ be a smooth, projective variety and $C\subset X$
an algebraic curve. It is not hard to see that there is
always a  meromorphic  volume form $\omega_m$ that is defined and nonzero
at all but finitely many points of $C$. We define the
{\it degree} of  $\omega_m$ on $C$ as
$$
\deg_C\omega_m:=
\#(\mbox{zeros of $\omega_m$ on $C$})-\#(\mbox{poles of $\omega_m$ on $C$}),
$$
where both zeros and poles are counted with multiplicities.
\end{defn}

The Chern form and the algebraic degree are connected by the 
Gauss--Bonnet theorem.

\begin{thm} \label{gauss-b.thm}
Let $X$ be a smooth, projective variety.
Let $\omega_r$ be a $C^{\infty}$ complex volume form and
$\omega_m$ a meromorphic  volume form. Then, for every
algebraic curve $C\subset X$
$$
\int_C  c_1(X)=\int_C \tilde c_1(X, \omega_r)=-\deg_C\omega_m.
\eqno{(\ref{gauss-b.thm}.1)}
$$
\end{thm}

(The minus sign comes from the happenstance that differential geometers
prefer to work with the tangent bundle while the volume forms
use the (determinant of the) cotangent bundle.)

\subsection*{Positivity/negativity and complex differential  geometry}{\ }

In differential geometry it is especially nice to
work with metrics whose curvature is everywhere
positive (or everywhere zero or everywhere negative) but
these rarely exist. A usual weakening is to work with
K\"ahler metrics that satisfy the
{\it Einstein condition}: the Ricci curvature
should be a constant multiple of the metric;
see \cite[Chap.19]{MR2325093} for definitions and an introduction.

If this {\it Einstein constant} is positive, then
in (\ref{gauss-b.thm}.1) we integrate an everywhere
positive form. Thus $\int_C  c_1(X) $ is positive for every curve $C$.
We hope that in this case there are  meromorphic  volume forms
with  poles  (but no zeros). 

Similarly, if the Einstein constant is negative, then
in (\ref{gauss-b.thm}.1) we integrate an everywhere
negative form. Thus $\int_C  c_1(X) $ is negative for every curve $C$.
We hope that in this case there are  holomorphic  volume forms
 (usually with zeros).

Algebraic geometry can be used to understand the numbers
$\deg_C\omega_m$, hence the values of the integrals  $\int_C  c_1(X) $.
It is a very difficult task to use the
positivity/negativity of the  integrals  $\int_C  c_1(X) $
to obtain a  K\"ahler  metric with  positive/negative Einstein constant.

For smooth varieties Aubin and Yau  proved existence in 1977
when  $\int_C  c_1(X) $ is always negative or when
 $\int_C  c_1(X) $ is always zero.
The singular case is treated in \cite{MR2505296, ber-g}. 
The positive curvature case  is more
subtle; a complete answer is not yet known.

While our approach to the structure of varieties is guided by
these curvature considerations, in algebraic geometry we can understand
only the algebraic degree of the first Chern class. Thus we look at the
functional
$$
C\mapsto \int_Cc_1(X)  
$$
and focus on those varieties where this is  everywhere
negative (or everywhere zero or everywhere positive).

The Main Conjecture then asserts that every variety can be built up
from these special varieties in a rather clear process.

\section{The Main Conjecture}\label{blocks.sec}
%The classification of algebraic varieties

On a typical variety $X$, the Chern class $c_1(X)$ is positive on
some curves and negative on others, in a rather unpredictable way.
Using the first Chern class and the theory of algebraic surfaces as our guide,
we focus on three basic ``especially simple'' types of 
smooth, projective  varieties.
These are the  ``building blocks''  of all algebraic varieties.

\begin{say}[{\bf Negatively curved}]\label{blocks.ssc.1}
 These are the varieties where
 $\int_Cc_1(X)$ is negative for every curve $C\subset X$.
This is the largest class of the three.
\end{say}

\begin{say}[{\bf Flat or Calabi--Yau}]\label{blocks.ssc.2}
 Here  $\int_Cc_1(X)$ is zero for every curve $C\subset X$.
They play an especially important role in string theory and mirror symmetry;
see \cite{MR1711184, MR2003030} for introductions.
\end{say}

\begin{say}[{\bf Positively curved or Fano}]\label{blocks.ssc.3}
 Here $\int_Cc_1(X)$ is positive for every curve. 
 There are  few of these varieties in each dimension, but
they occur especially frequently in applications.
\end{say}

A  simple set of examples to keep in mind is the following.
A smooth hypersurface  $X_d\subset \cc\pp^n$ of degree $d$ is 
negatively curved if  $d>n+1$, flat if $d=n+1$ and 
positively curved if $d<n+1$.

A variety in any of these 3 classes is considered ``simplest,''
but we do not yet have enough ``simplest'' varieties for 
answering the first Main Question. % \ref{prob1}.
For example, taking products of these we get examples where  $c_1(X)$ 
has different signs on different curves. Two of these
possible ``mixed types''  are relevant for us.

Consider a product  $X:=N\times F$ of a negatively curved and of a flat variety.
It is clear that $\int_Cc_1(X)\leq 0$ for every curve $C\subset X$
and $\int_Cc_1(X)= 0$ only if $C$ lies in a fiber of the first projection
$N\times F\to N$. This observation leads to the 4th class.

\begin{say}[{\bf Semi-negatively curved or Kodaira--Iitaka type}]
\label{blocks.ssc.4}
Here $\int_Cc_1(X)\leq 0$ for every curve $C\subset X$
and there is a  unique morphism  $I_X:X\to I(X)$ such that 
$\int_Cc_1(X)=0$ iff $C$ is contained in a fiber of  $I_X$.

This includes the classes \ref{blocks.ssc.1}--\ref{blocks.ssc.2}:
 $I_X$ is an isomorphism for negatively curved varieties and
a constant map in the flat case.

In the intermediate cases, when
$0<\dim I(X)<\dim X$, almost all fibers of $I_X$ are
 Calabi--Yau varieties. Thus one can view these as
families of (lower dimensional) Calabi--Yau varieties parametrized by  
the (lower dimensional) variety $I(X)$. If we understand families
of (lower dimensional) varieties well enough, we understand $X$.
(This is one of the reasons  we are interested in
the second Main Question.)
Furthermore, in these cases $I(X)$ is negatively curved
in a ``suitable sense,'' though we do not yet have a
final agreed-upon definition of what this  means. 
\end{say}

Next consider a product  $X:=N\times P$ of a negatively curved and of a 
positively curved variety. If a curve  $C$ 
lies in a fiber of the first projection then 
 $\int_Cc_1(X)> 0$, but there are many other such curves.
%Similarly, the fibers of the second projection can not be characterized
%by  $\int_Cc_1(X)< 0$. 
Nonetheless, the first projection 
is uniquely determined by $X$ and this 
 leads to the definition of the 5th class.

\begin{say}[{\bf Positive fiber type}]\label{blocks.ssc.5}
I really would like to say that in these cases
 there is a  unique morphism  $m_X:X\to M(X)$ such that 
$M(X)$ is semi-negatively curved and 
$c_1(X)$ is positive on all the fibers.
(To avoid trivial cases, we also assume that $\dim M(X)<\dim X$.)
This, unfortunately, still does not give  enough ``simplest'' varieties for 
the first Main Question. % \ref{prob1}. 
It took quite some time to arrive at the correct definition,
to be discussed  in Section \ref{ratconn.sec}.
\end{say}

We can now  state a precise version of
the first Main Question.

\begin{main-conj} \label{main-conj.1}
Every algebraic variety $X$ is birational to a variety
$X^{\rm m}$ that is either of type 
{\rm (\ref{blocks.ssc.4})} or  of type {\rm (\ref{blocks.ssc.5})}.
\end{main-conj}

\noindent{\bf Complement.}   $X^{\rm m}$ 
-- especially in case
(\ref{blocks.ssc.4}) -- is called a
{\it minimal model} of $X$. 

 In the semi-negatively curved case  $I\bigl(X^{\rm m}\bigr)$
is unique but $X^{\rm m} $ itself is not. However, it is quite well understood
how the different  $X^{\rm m} $ are related to each other.
(This is the story of {\it flops,} see \cite{MR1159257, MR2827807}.)
By contrast, in   case (\ref{blocks.ssc.5})
 it is very hard to determine
when two such varieties $X^{\rm m}_1 $ and $X^{\rm m}_2 $ are birational.
\smallskip

\noindent{\bf Caveat.}  While the Main Conjecture  is expected to be true, 
in general
one has to allow {\it terminal} singularities --
 to be defined in (\ref{terminal.defn}) --
 on $X^{\rm m}$. 

This was a rather
difficult point historically 
since over a century of experience suggested
that singularities should be avoided. 
For surfaces  terminal = smooth, thus the issue of singularities
did not come up in Theorem \ref{min.mod.surf.thm}.

By now the correct classes of
singularities have been established and, for many questions we consider,
 they do not seem to cause any problems. We describe these singularities
in  Section \ref{sing.sec}.

\begin{say}[Traditional names] \label{trad.names.say}
A variety $X$ is said to be
of {\it general type} if $\dim I\bigl(X^{\rm m}\bigr)=\dim X$.
In this case $X\map I\bigl(X^{\rm m}\bigr)$ is birational and 
$I(X):=I\bigl(X^{\rm m}\bigr)$ is  called the
{\it canonical model} of $X$;  it 
has canonical singularities (\ref{can-lc.defn}).
We see in Section \ref{moduli.sec}
 that   the second Main Question %\ref{prob2}
 has a good answer for families of 
canonical models.

The {\it Kodaira dimension} of a variety $X$ is the dimension
of $I\bigl(X^{\rm m}\bigr)$. 

The Kodaira dimension is defined to be
 $-\infty$ for the  class (\ref{blocks.ssc.5}).

The Main Conjecture is usually broken down into two parts
that are, in principle, independent of each other.
The first part  separates the classes \ref{blocks.ssc.4}\ and \ref{blocks.ssc.5}
from each other  and the second part
provides the structural description in case \ref{blocks.ssc.4}.
These forms first appear in Reid's paper \cite[Sec.4]{r-mmc3}.

\medskip
\noindent \ref{trad.names.say}.1 
{\it Minimal Model Conjecture.}  % \label{mmp.conj}
Every algebraic variety $X$ is birational to a
 variety
$X^{\rm m}$ such that
 either $c_1\bigl(X^{\rm m}\bigr)$ is semi-negative
 or there is a morphism  to a lower dimensional variety 
 $\pi:X^{\rm m}\to S$ such that
$\int_Cc_1\bigl(X^{\rm m}\bigr)>0$ if $C$ is contained in a fiber of  $\pi$.
(In the second case the map $\pi$  need not be 
unique and it does not give the best structural description.)

\medskip
\noindent \ref{trad.names.say}.2
 {\it Abundance Conjecture.}   %\label{abund.conj}
If   $c_1(Y)$ is semi-negative then
 there is a  unique morphism  $I_Y:Y\to I(Y)$  such that 
$\int_Cc_1(Y)=0$ iff $C$ is contained in a fiber of  $I_Y$.
\end{say}

\section{Rationally connected varieties}\label{ratconn.sec}

Before we consider  minimal models, we describe
the structure we expect for  varieties in the 
5th class (\ref{blocks.ssc.5}).
An introduction aimed at non-specialists is given in
\cite{k-simplest}. More detailed accounts are in
\cite{ar-ko, rc-book}.

Clebsch and Max Noether  noticed around 1860--1870
 that, when the numerical invariants suggest that 
a surface could be birational to $\cc\pp^2$, then it is.
The final result along these lines was established by
Castelnuovo in 1896.

Analogous questions in higher dimension turned out to be
much harder. Fano classified smooth positively curved 3--folds
around 1930. (He missed some cases though, so did subsequent ``complete''
lists produced in the 1970's and then in the 1980's. The (hopefully)
final list was not established until 2003.) 
This is, however, one area where the singularities do matter;
we still do not know  all positively curved 3--folds
with terminal singularities. %; see \cite{MR1668579} for a survey.

It appears that instead of global descriptions
we should focus on  {\it rational curves} in a variety;
these are the images of  morphisms  $\phi:\cc\pp^1\to X$.
For a projective variety $X$, the following dichotomy is quite easy
to establish.

i) either the rational curves cover a  subset of $X$
which is {\it meager} 
 (that is, a 
 countable union of nowhere dense closed subsets)

ii) or the rational curves cover all of $X$.

These two cases correspond to the alternatives in
the  Main Conjecture. % \ref{main-conj.1}.
That is, if
 $X$ is  birational to a semi-negatively curved variety
 then rational curves cover a meager subset and, conjecturally,
the converse also holds.

The correct approach to  the best structural description
of the 5th class \ref{blocks.ssc.5}
was  not discovered until 1992
(Koll\'ar--Miyaoka--Mori \cite{KMM92a}).
The key observation is that
we should even change the class \ref{blocks.ssc.3}.
 Instead of a curvature description,
we should focus  on  rational curves contained in a variety.

\begin{defn} 
A  projective variety $X$ is called 
{\it rationally connected} if, for any number of points  $x_1,\dots, x_r\in X$,
 there is a
morphism $\phi:\cc\pp^1\to X$ whose image passes through $x_1,\dots, x_r$.
\end{defn}

I claim that rationally connected varieties constitute the
``correct'' birational version of being
positively curved. This is not a precise mathematical assertion
since not every rationally connected variety is birational to
a positively curved variety, not even when singularities are allowed.
Rather, the assertion is that  any answer to  the first Main Question 
 needs to work with rational connectedness
instead of positivity of curvature.

\begin{say}[Supporting evidence]\label{supp.evid.rc.say}{\ }

  It is easy to see that $\cc\pp^n$ is rationally connected.
More generally, every positively curved variety is 
rationally connected  (Nadel \cite{Nadel91}, Campana \cite{Campana92b}, 
Koll\'ar--Miyaoka--Mori \cite{KMM92c}, Zhang \cite{MR2208131}). 
%\item Being rationally connected is a birational invariant.

 Being rationally connected is invariant under smooth deformations
and birational maps \cite{KMM92c}.

 Rationally connected varieties share key arithmetic properties of
rational varieties over $p$-adic fields (Koll\'ar \cite{k-loc}),
  finite fields (Koll\'ar--Szab\'o \cite{MR2019976}, Esnault 
\cite{MR1943746})
and  function fields of curves  (Graber--Harris--Starr \cite{ghs},
 de~Jong--Starr \cite{MR1981034}).

 The loop space of a rationally connected variety is
also rationally connected  (Lempert--Szab\'o \cite{MR2372727}).
\end{say}

The notion of rational connectedness allows us to give the
correct description of the  class \ref{blocks.ssc.5}.
A weaker variant  is proved  in \cite{KMM92a};
the form below combines this with  \cite{ghs}. 

\begin{thm} Let $X$ be a  variety that is covered by rational curves.
Then  there is a unique
(up to birational equivalence) map
$m_X:X\map M(X)$ such that
\begin{enumerate}
\item almost all fibers of $m_X$ are rationally connected and
\item rational curves cover only a meager subset of $M(X)$. 
\end {enumerate}
\end{thm}

There are 
two main open geometric problems about rationally connected varieties.
The first concerns a topological characterization.
In its naive form the question asks: What can we tell about a variety
from its underlying topological space? It seems that the answer is: not much.
However, the underlying topological space of a smooth variety
carries  a natural symplectic structure and this seems to incorporate
much more information.

\begin{conj} {\rm \cite[Conj.4.2.7]{k-lowdeg}} Being rationally connected is a 
property of the
underlying symplectic structure.
\end{conj}

For partial results see \cite{k-lowdeg, MR2521651}.

The other problem asks if we could strengthen the
definition of rationally connected varieties.
Note that $\cc\pp^n$ contains not just many rational curves
but also many higher dimensional rational 
subvarieties  (hyperplanes, hyperquadrics, ...). 
Maybe this is also a general property of rationally connected varieties?
As far as I know, 3-dimensional rationally connected varieties
always contain rational surfaces. I believe, however,
 that this is not the case in higher dimension.

\begin{conj} {\rm \cite[Prob.56]{k-simplest}} 
Many rationally connected varieties
do not contain any rational surface.
\end{conj}

\section{Minimal Models}

This is a short history of {\it Mori's program}, also called   
the {\it Minimal Model Program} and frequently abbreviated as {\it MMP}.
For general introductions see \cite{CKM88, km-book}
or the technically more detailed \cite{ka-ma-ma, MR2359340, hac-kov}. 

\begin{say}[{\bf Iitaka's program, 1970--85}]
This approach predates
the Main Conjecture. At the beginning it was not
even suspected that the Main Conjecture could be true, in fact,
lacking the right class of singularities,
it was assumed that the Main Conjecture
 would fail for most varieties. Thus the aim of Iitaka's program was
to sort varieties into 5 broad types that (as we now  know) exactly
correspond to the ones  in (\ref{blocks.ssc.1}--\ref{blocks.ssc.5}).
 The main contributors were,
in rough historical order, 
Iitaka, Ueno, Fujita, Kawamata, Viehweg and Koll\'ar;
see \cite{ueno, Mori87} for surveys.
\end{say}

\begin{say}[{\bf Canonical and terminal singularities, Reid 1980--83}]
 Reid was studying higher dimensional
analogs of Du~Val singularities of surfaces (\ref{duval.say});
 obtaining rather complete
descriptions in dimension 3.   %\cite{r-c3f, r-mmc3}.  
It was quite important that when
Mori's program lead to singularities, the relevant classes were
already there and were known to be well behaved.
An especially readable account is 
\cite{r-ypg}.
\end{say}

\begin{say}[{\bf The birth of Mori's program, 1981--88}]

Mori's groundbreaking paper \cite{Mori82} introduces 3 new ideas.

 If $c_1(X)$ is not semi-negative then, by definition, 
$c_1(X)$ is positive on some curve $C\subset X$.
Mori first  proves that there is such a rational curve;
that is, there is a  morphism  $\phi:\cc\pp^1\to X$ such that
$c_1(X)$ is positive on its image.
 It is quite remarkable that
the proof goes through algebraic geometry over finite fields.
To this day there is no proof known that avoids this; in particular
this step is not yet known for complex manifolds that are not algebraic.

Second, he identifies the ``most positive'' such maps 
 $\phi:\cc\pp^1\to X$; this is called {\it extremal ray} theory.

Third, in dimension 3 he gives a complete description of
all extremal rays and the resulting  map $X\to X_1$
that removes the ``most positive'' part of $X$.

The program now seems clear (at least in dimension 3).
 Repeat the procedure for $X_1$
and prove that after finitely many steps we end up
with   $X\to X_1\to\cdots\to X_r$ such that $c_1(X_r)$ is  semi-negative.
This is called {\it Mori's program} or
 {\it Minimal Model  program.}

There are two, rather formidable, problems.
In many cases the new variety $X_1$ is smooth but 
sometimes it is singular. Luckily, these singularities have been
studied by Reid, at least in dimension 3.
Still, it is necessary to
establish the above 3 steps for singular varieties.
This was accomplished rather rapidly by 
Kawamata, Reid,  Shokurov and  Koll\'ar.
The program was first  written down in  \cite[Sec.4]{r-mmc3}.

The more serious problem is that in some cases
taking the contraction  $X_i\to X_{i+1}$ is clearly not the right step.
Instead we have to take a step back and construct a new
variety   $X_i^+ $ that sits in a {\it flip diagram}
$$
\xymatrix{%
X_i  \ar@{-->}[rr]^{\phi_i} \ar[rd]_(.45){p_i}
&& X_i^+  \ar[ld]^(.35){p_i^+} \\
& X_{i+1} &
}
$$
%$$
%\begin{array}{rcl} X_i &\stackrel{\phi_i}{\map} &X_i^+\\ 
%p_i &\searrow \quad\swarrow & p_i^+\\  &X_{i+1}&
%\end{array}
%\eqno{(\ref{mmp.restricts.do.divs.1}.1)}
%$$
Geometrically, we start with $X_i$, find an especially
badly behaving  $\cc\pp^1\cong C_i\subset X_i$ and remove it. 
Then we compactify the resulting $X_i\setminus C_i$
by attaching  another curve  $C_i^+\cong \cc\pp^1$ but differently.
The key difference is a sign change:
$$
\int_{C_i} c_1\bigl(X_i\bigr)>0\qtq{but} \int_{C_i^+} c_1\bigl(X_i^+\bigr)<0.
$$
This operation is called a {\it flip.}
For more about flips, see \cite{MR1159257, MR2827807}.

Flips are reminiscent of  {\it Dehn surgery}  in  3--manifold topology
where we remove a circle and put it back differently.

In dimension 3 the existence of flips is proved in a very difficult paper
by Mori \cite{mori-MR924704}, which completes the program in this case.
A detailed description of 3--dimensional flips is given in
\cite{MR1149195}. The list is rather lengthy; this makes it unlikely that
a similarly complete answer will ever be worked out in higher dimensions.
\end{say}

\begin{say}[{\bf Log variants: Kawamata,  Shokurov, 1984--1992}]
The Iitaka program established that for many results one can
work with cohomology classes in $H^2(X,\rr)$ that are close enough to
 the first Chern class.
This turned out to be a very powerful tool. By choosing
the perturbations appropriately, we can focus our attention
on one or another part of a variety.  These are somewhat technical questions
but by now we understand how to work with them and
most  applications of the Minimal Model Program use a perturbed case. 
\end{say}

\begin{say}[{\bf Abundance: Kawamata, Miyaoka, 1987--1992}]
Even for surfaces, the Abundance Conjecture  \ref{trad.names.say}.2
is a rather subtle result.
It is even harder for 3--folds. The proofs use many special
properties of surfaces; this is why the higher dimensional cases
are still not well understood.
A rather complete account of the 3--dimensional 
methods is given in \cite{k-etal}.
\end{say}

\begin{say}[{\bf Inductive approach in low dimensions: Shokurov, 1992--2003}]
 In retrospect, the
key development of the decade was an inductive approach to flips.
A detailed treatment of the 3--dimensional case is given in
\cite{k-etal}. For the rest of the nineties progress was slow, culminating
in a treatment of 4--dimensional flips. 
There were many technical difficulties to overcome and the
importance of these methods   was not fully appreciated at
first since the dimension reduction leads to a much more complicated
problem that seems to fail in higher dimensions.
\end{say}

\begin{say}[{\bf The Corti seminar, 2003--2005}]
Over the course of several years a group led by Corti
developed the previous ideas further and integrated them with the
rest of the program \cite{MR2359340}.
 This provided the bridge to the general case.
\end{say}

\begin{say}[{\bf The general type case: Hacon and M\textsuperscript{c}Kernan, 
 2005--2010}]
The real  breakthrough was achieved in \cite{MR2359343} 
where the existence of flips in dimension $n$ was reduced to
an instance of the MMP in dimension $n-1$. This left a
series of global questions to resolve. The paper
\cite{bchm} settled everything for varieties of general type.
A good introduction is in  \cite{MR2743824}.

At about the same time 
Siu started to develop an analytic approach which aims to get
$I(X^{\rm m})$, without going through the
individual steps; see
\cite{MR2827808} for an overview.
An algebraic variant of this is in \cite{cas-laz}.

\end{say}

\begin{say}[{\bf Abundance: Hacon and Xu, 2012--}]
Although the Abundance conjecture  is known in very few cases,
there has been significant progress when $\dim I(X)$ is 
expected to be close to $\dim X$. The log version of the  special case when 
$\dim I(X)=\dim X$ is especially important for applications
in moduli theory. These have been settled in \cite{MR3032329, MR2955764}.
\end{say}

\begin{say}[{\bf Positive characteristic,  Hacon and Xu, Birkar,  Patakfalvi,
2012--}]
Mori's original works are very geometric
and these ideas quickly lead to a simple proof of the 
2-dimensional case of the Main Conjecture  in positive characteristic.
However,  subsequent developments
rely very heavily on Kodaira-type vanishing theorems that
are known to fail in positive characteristic, although no actual
failure is known in the cases used by the program.
The 3--dimensional case was recently settled in
\cite{2013arXiv1302.0298H, 2013arXiv1311.3098B}.
Substantial parts of the Iitaka program are proved
 in positive characteristic in  \cite{2013arXiv1308.5371P}.
\end{say}

\begin{say}[{\bf Open problems}] From our point of view, the main open problem
is to complete the missing parts of the Main Conjecture.

It is known that the MMP always runs, 
that is, the sequence of contractions and flips
$X=X_1\map X_2\map \cdots$ exists. 
The problem is that
it is not clear how to prove that the process eventually stops.
In the 3--dimensional case, Mori's approach provides a
rather complete description of the steps of the MMP.
This gives many ways to show that each step
improves various invariants and that eventually the process stops.
By contrast, the method of Hacon--M\textsuperscript{c}Kernan produces the steps
of the MMP in a rather indirect way. We have very little
information about the steps beyond their existence.
\end{say}

\section{Singularities of the Minimal Model Program} \label{sing.sec}

So far we have been sweeping the singularities of the
minimal models under a rug, but it is time for a  look at them.
Understanding the correct class of singularities
is crucial in the development of the structure theory of algebraic varieties.
This is a somewhat technical subject with many
difficult questions and methods but by now we understand
these singularities well enough that in many questions 
they do not cause any problems. A rather complete treatment is given in
\cite{kk-singbook}. Here I focus on the main ideas behind the definitions.

Given a variety $Y$, one frequently looks at a resolution
of singularities  $f:X\to Y$ as in Theorem \ref{res.sing.thm}
and translates problems on $Y$ to
questions on $X$. Then the hard part is to interpret the answer obtained on
$X$ in terms of $Y$. Here the key seems to be the inverse function
theorem.

\begin{say}[{\bf The inverse function theorem}]
The classical inverse function theorem says that if
${\mathbf f}:=(f_1,\dots, f_n):\rr^n_{\mathbf x}\to \rr^n_{\mathbf y}$
is a differentiable map then ${\mathbf f} $
has a local inverse at a point $p\in \rr^n_{\mathbf x}$
 iff the Jacobian determinant
$$
\operatorname{Jac}({\mathbf f}):=
\det\left(\frac{\partial f_i}{\partial x_j}\right)
$$
does not vanish at $p$. We can also think about it 
in terms of the ``standard'' volume forms
$\omega_{\mathbf x}:=dx_1\wedge\cdots\wedge dx_n$ and
$\omega_{\mathbf y}:=dy_1\wedge\cdots\wedge dy_n$.
Then
$$
{\mathbf f}^*\omega_{\mathbf y}=
\operatorname{Jac}({\mathbf f})\cdot \omega_{\mathbf x},
$$
thus the vanishing/non-vanishing of the Jacobian tells us
how the pull-back of the ``standard'' volume form of the target
compares to the ``standard'' volume form of the source.

Note that the Jacobian itself depends on the choice of the coordinates,
but its vanishing or non-vanishing depends only on ${\mathbf f} $. 

In the complex analytic setting
one can use  the ``standard'' complex volume forms
$\omega_{\mathbf z}:=dz_1\wedge\cdots\wedge dz_n$ and
$\omega_{\mathbf w}:=dw_1\wedge\cdots\wedge dw_n$
on the unit balls  $\bb^n_{\mathbf z}\subset \cc^n_{\mathbf z}$
and $\bb^n_{\mathbf w}\subset \cc^n_{\mathbf w}$.
Given a holomorphic map
${\mathbf f}:=(f_1,\dots, f_n):\bb^n_{\mathbf z}\to \bb^n_{\mathbf w}$
we get  that
$$
{\mathbf f}^*\omega_{\mathbf w}
=\det\left(\frac{\partial f_i}{\partial z_j}\right)\cdot \omega_{\mathbf z}
=:\operatorname{Jac}({\mathbf f})\cdot \omega_{\mathbf z},
$$
and ${\mathbf f} $
has a local inverse 
 iff $\operatorname{Jac}({\mathbf f}) $
does not vanish at $p$. 
\end{say}

\begin{say}[{\bf The Jacobian in the singular case}]\label{jac.sing}

Let $X$ be a normal  algebraic variety 
and $p\in X$ a singular point. It is quite easy to see
that if $\omega_1, \omega_2$ are two holomorphic volume forms on 
$X\setminus\sing X$ in  a neighborhood of a singular point $p\in \sing X$
 then
there is a unique holomorphic function $\phi$ such that 
$\omega_1=\phi\cdot\omega_2$ and $\phi(p)\neq 0$.
Thus all holomorphic volume forms on $X\setminus\sing X$
 have the same asymptotic
behavior near   $\sing X$. 
The local existence of such forms is a slightly technical question, 
so let us just focus on an example.
If
$Y=\bigl(f(w_1,\dots, w_{n+1})=0\bigr)\subset \cc^{n+1}$
is a hypersurface then the ``standard'' volume form is given by
$$
\omega_Y=(-1)^i
\frac{dw_1\wedge \cdots\wedge dw_{i-1}\wedge dw_{i+1}\wedge\cdots\wedge dw_{n+1}}
{\partial f/\partial w_i}.
\eqno{(\ref{jac.sing}.1)}
$$
(It is easy to check that this is independent of $i$. Note also
that $\omega_Y $ is not defined when all of the
$\partial f/\partial w_i $ vanish; which happens  exactly on $\sing Y$.) 
Thus if $f:\bb^n_{\mathbf z}\to Y$ is 
holomorphic then we can define the {\it Jacobian} of $f$
by the formula
$$
\operatorname{Jac}(f):=\frac{f^*\omega_Y}{ \omega_{\mathbf z}}.
$$
Note that due to the denominators in (\ref{jac.sing}.1), 
in general $\operatorname{Jac}(f) $ can have poles.

For example, consider the singularity
$Y_{n,d}:=\bigl(w_1^d+\cdots+w_n^d=w_{n+1}^d\bigr)\subset \cc^{n+1}$
and a holomorphic map $f:\bb^n_{\mathbf z}\to Y$ given by
$$
f:(z_1,\dots, z_n)\to 
\bigl(z_1, z_1z_2, \dots, z_1z_n, z_1\sqrt[d]{1+  z_2^d+\cdots+z_n^d}\bigr).
\eqno{(\ref{jac.sing}.2)}
$$
Then $\omega_{Y_{d,n}}=-d^{-1}w_1^{1-d} dw_2\wedge \cdots\wedge dw_{n+1}$ and
we easily compute  that  the Jacobian of $f$ has a zero/pole of order  $n-d$
along the hyperplane $(z_1=0)$.

As in the classical case, the Jacobian of $f$ depends on the choice of
the ``standard'' volume forms but the vanishing/non-vanishing
or the order of vanishing of the Jacobian depends only on $f$.

We can now define terminal singularities; these form the smallest
possible class needed for the Main Conjecture. %  \ref{main-conj.1}.
\end{say}

\begin{defn}\label{terminal.defn}
 A normal variety  $Y$ has
{\it terminal} singularities iff the 
inverse function theorem holds for $Y$. That is, if 
$f:\bb^n_{\mathbf z}\to Y$ does not have a local inverse at $p\in \bb^n_{\mathbf z}$
then  $\operatorname{Jac}(f)$ vanishes at  $p$.
(There is a small  problem when the exceptional set of $f$ is
too small, we can  ignore it for now.)
\end{defn}

For canonical models and for moduli questions, two more types of
singularities are needed.

\begin{defn}\label{can-lc.defn}
 A normal variety  $Y$ has
 {\it canonical} singularities iff 
 $\operatorname{Jac}(f) $ is holomorphic for every $f:\bb^n_{\mathbf z}\to Y$
and  {\it log-canonical} singularities iff 
 $\operatorname{Jac}(f) $ has at most simple poles
for every $f$.
\end{defn}

The above computations suggest  (and it is indeed true) that
  $Y_{n,d}$ (as in \ref{jac.sing}.2) is terminal iff $d<n$, canonical iff
$d\leq n$ and log  canonical iff
$d\leq n+1$.

\begin{say}[Local volume of $Y$ near $\sing Y$] A good way to think about
these singularities is as follows. Pick a point $p\in \sing Y$
and let $\omega_Y$ be a ``standard'' local complex volume form.
Then $(\sqrt{-1}/2)^n \omega_Y\wedge\bar\omega_Y$ is a real  volume form
and we can ask about the
{\it local volume of $X$}, that is, 
$\int_U (\sqrt{-1}/2)^n \omega_Y\wedge\bar\omega_Y$
for a suitably small  neighborhood $p\in U\subset X$.

If $Y$ has a canonical singularity near $p$ then the
local volume is finite. In the log-canonical case
the local volume is infinite but barely. If $g$ is any
holomorphic function vanishing on $\sing Y$ then
$\int_U|g|^{\epsilon}(\sqrt{-1}/2)^n \omega_Y\wedge\bar\omega_Y$
is finite for every $\epsilon>0$.
\end{say}

\begin{say}[{\bf Intermediate differential forms}]

On an $n$-dimensional variety we have so far considered
holomorphic $n$-forms only but for several questions
one also needs to understand the pull-back $f^*\eta$ of 
lower degree differential forms
as well. This proved to be surprisingly difficult but
almost all local questions were settled by 
Greb--Kebekus--Kov\'acs--Peternell \cite{GKKP10}. 
\end{say}

\section{Moduli of varieties of general type}\label{moduli.sec}

Let ${\mathbf X}$ be a class of projective  varieties,
for instance curves or surfaces of a certain type.
The theory of moduli aims to  find ``optimal'' ways
to write down all varieties in the class  ${\mathbf X}$.

This is a  large theory with many aspects. The 3 volumes of
\cite{mod-hand} contain surveys of most of the active areas.
Here my aim is to focus on just one of them: the
moduli of varieties of general type. 
Introductions are given in 
\cite{MR2483953, hac-kov, k-modsurv} while a detailed treatment should be
in \cite{k-modbook}.

We start with the historically first example.

\begin{exmp}[Elliptic curves]  They can all be given by
an affine equation
$$
E(a,b,c):=(y^2=x^3+ax^2+bx+c)\subset \cc^2;
$$
the corresponding projective curve has a unique 
 point $[p]$  at infinity.
Here $c_1(E)=0$, so it is best to think of this as
elliptic curves with a marked point  $[p]$.
The curve $E(a,b,c) $ is smooth iff the 
discriminant of the cubic
$$
\Delta(a,b,c):=18abc-4a^3c+a^2b^2-4b^3-27c^2\qtq{is not zero.}
$$
Two such curves  %$E(a_1,b_1,c_1) $ and $E(a_2,b_2,c_2) $
are isomorphic iff there is an affine-linear transformation
$(x,y)\mapsto (\alpha^2 x+\beta, \alpha^3 y)$
 that transforms one equation into the other.
All these transformations form a (2-dimensional)   group $G$. 
Thus we get the following.

{\it Version 1.} 
The isomorphism classes of all elliptic curves are in one-to-one
correspondence with the orbits of $G$ on
$\cc^3\setminus(\Delta(a,b,c)=0)$.

Next we need to identify the $G$-orbits.
The key is the  $j$-invariant 
$j\bigl(E(a,b,c)\bigr):= 2^{8}(a^2-3b)^3/
\Delta(a,b,c)$.
(The factor  $2^{8}$ is not important for us, it is 
there for number-theoretic reasons.) 
It is  not very hard to work out the following.

{\it Version 2.} 
Two  elliptic curves are isomorphic iff they have the same 
$j$-invariant. 

We can restate this as follows:

{\it Version 3.} 
The moduli space of elliptic curves is the complex line 
${\mathcal M}_1\cong \cc$
and the value $j(E)$ of the $j$-invariant gives the point
in ${\mathcal M}_1$ that corresponds to $E$.

The only sensible compactification of $\cc$ is $\cc\pp^1$,
so what corresponds to the point at infinity? This should be
a curve where the discriminant  of the cubic $x^3+ax^2+bx+c$  vanishes.
That is, when $x^3+ax^2+bx+c $ has a multiple root.
There are 2 types of such cubics. If there is a triple root
we get  $y^2=x^3$, a cuspidal curve.  If there is a double root
we get  $y^2=x^3+x^2$, a nodal curve. In this case the correct choice is to
go with the nodal curve.
\end{exmp}

\begin{say}[The main steps of a moduli theory] \label{main.steps.moduli}
We hope to do something similar with more general
algebraic varieties. We proceed in several steps.

{\it Step 1.} Identify a class of projective  varieties ${\mathbf X}$
that should have a ``good''  moduli theory. 
We aim to prove that such a theory exists for
negatively curved varieties as in (\ref{blocks.ssc.1}).
We allow canonical singularities, thus this includes
canonical models of varieties of general type.
(It seems that in most other 
cases there is no  ``good''  compactified moduli theory, unless
some additional structure is added on, for instance an ample divisor
as in \cite{ale-abvar}.)

{\it Step 2.} Add some extra data (also  called rigidification) first.
A typical extra datum 
is an embedding
$j:X\into \pp^N$ for some $N$.
 Use the additional data to get a  moduli space 
with a universal family 
$$
{\mathbf U}_{\mathbf X,j}\subset \pp^N\times  {\mathbf M}_{\mathbf X,j}
\qtq{with projection}
\pi_{\mathbf X,j}:{\mathbf U}_{\mathbf X,j}\to {\mathbf M}_{\mathbf X,j}
$$
such that
 every pair  $(X,j)$ occurs exactly once among the fibers of $\pi_{\mathbf X,j}$.
(It is not easy to show that one can choose a fixed $N$ that works for all
varieties in a given class. For smooth varieties 
this was proved by Matsusaka in 1972; the general case  was settled
recently by Hacon and Xu.)

{\it Step 3.} Next we  get rid of the extra data. Usually
 we have to take a quotient
by a Lie group  like  $\GL(N+1, \cc)$. 
This can be  hard but, if everything works out, 
at the end we have
$$
{\mathbf U}_{\mathbf X}:={\mathbf U}_{\mathbf X,j}/\GL(N+1, \cc),\quad
{\mathbf M}_{\mathbf X}:={\mathbf M}_{\mathbf X,j}/\GL(N+1, \cc)
$$
and a morphism
$\pi_{\mathbf X}:{\mathbf U}_{\mathbf X}
\to {\mathbf M}_{\mathbf X}$. (See Step 6 for the possible dependence on $N$.)

{\it Step 4.} In almost all cases, the resulting spaces are not compact
and compactifying them in a ``good'' way is difficult.
The key step is to
identify the limits of families of   varieties in ${\mathbf X}$
that should give a ``good'' compact moduli theory.
There is no a priori reason to believe that such a choice exists
or that it is unique. Finding the right choice in higher dimension
was the last conceptual step in the program.
For  canonical models of varieties of general type we  have the
``right'' answer, see (\ref{mod.KSB}) and (\ref{mod.gen}).

{\it Step 5.} We have to go back and redo Steps 1--3 for this more
general class of objects to get a compactified moduli theory
$$
\bar\pi_{\mathbf X}:\overline{\mathbf U}_{\mathbf X}
\to \overline{\mathbf M}_{\mathbf X}.
$$

{\it Step 6.} An extra issue that arises is that
 the compactifications could also  depend on the 
dimension of the $\pp^N$ chosen in Step 2. This does not seem to happen
for ${\mathbf M}_{\mathbf X}$ itself (at least for $N$ large enough)
but it does happen for $\overline{\mathbf M}_{\mathbf X} $
for some of the proposed
variants.

{\it Step 7.} Finally, if everything works out,
we would like to study the properties of ${\mathbf M}_{\mathbf X}$,
 $\overline{\mathbf M}_{\mathbf X}$ and to use these to prove further theorems.
\end{say}

Next we review the historical development of the
higher dimensional theory.

\begin{say}[{\bf Geometric Invariant Theory: Mumford, 1965}]
Riemann probably knew that all smooth, compact Riemann surfaces of
a given genus $g$ form a nice family, but
the  moduli spaces ${\mathcal M}_g$ were first rigorously constructed by
Teichm\"uller in 1940 as an analytic space and by
Mumford in 1965 as an algebraic variety.
Mumford's book \cite{git} presents a program to construct moduli spaces
under rather general conditions and uses it to obtain ${\mathcal M}_g$.
Using these methods, moduli spaces were
constructed for surfaces (Gieseker, 1977) and  
for higher dimensions (Viehweg, 1990).

The correct compactification of these  moduli spaces
 was much less clear. In principle,
GIT provides an answer, but the resulting compactification 
might  depend on  the embedding dimension
chosen in (\ref{main.steps.moduli}.Step 2). 
Recently Wang--Xu \cite{xu-wang}
prove that, for surfaces and in higher dimensions, 
the GIT compactification  does depend on  the embedding dimension. 
(The current examples, however, do not exclude the possibility that
some variant of the GIT approach does provide an answer
that is independent of the embedding dimension.)
\end{say}

\begin{say}[{\bf Compact moduli of curves: Deligne and Mumford, 1969}]
\label{comp.dm.say}

The optimal compactification  of ${\mathcal M}_g$ is constructed
in \cite{del-mum}. In the boundary  
$\overline{\mathcal M}_g\setminus {\mathcal M}_g$
we should allow reducible curves  $C=\cup_i C_i$ that satisfy two
restrictions.

({\it Local property}) $C$ has only {\it nodes} as singularities.
In suitable local analytic coordinates these are given
by an equation $(xy=0)\subset \cc^2$.
As in  (\ref{jac.sing}) the ``standard'' volume form on a node is given by
$\frac{dx}{x}$  (on the line $(y=0)$) and by
 $-\frac{dy}{y}$  (on the line $(x=0)$). These forms
have a simple pole at the singularity, corresponding to the restriction on
log canonical singularities in (\ref{can-lc.defn}).

({\it Global property}) Instead of  each $c_1(C_i)$ being negative,
we assume that  each $c_1(C_i)-D_i$ is negative where $D_i$ is the sum of the 
nodes that lie on $C_i$. (Thus we allow $C_i\cong \cc\pp^1$, as long as at least
 3 nodes also lie on $C_i$.)
\end{say}

\begin{say}[{\bf Compact moduli of surfaces: Koll\'ar and Shepherd--Barron, 1988}]
\label{mod.KSB}
It was clear from the Mumford--Gieseker approach that one should work
with the {\it canonical models}    of  surfaces
of general type (as in \ref{duval.say}) in order to get a good moduli theory,
but the correct class of singular limits was not known.

An approach using minimal models was proposed in
\cite{ksb}: given a family of canonical models
over a punctured disc $S^*\to \Delta^*$, 
first construct any compactification whose
central fiber is a reduced simple normal crossing divisor
and then take the (relative) canonical model. It is not hard to see that
this gives a unique limit. This says that at the boundary 
of the moduli space
we should allow {\it stable surfaces:} 
 reducible surfaces  $S=\cup_i S_i$ that satisfy two
restrictions.

({\it Local property}) $S$ has so--called  {\it semi-log canonical}
singularities. What are these? First of all, aside from
finitely many points $S$ is either smooth or has two local branches meeting
transversally, like  $(xy=0)\subset \cc^3$. These are the natural
generalizations of nodes. Then we can have log canonical singularities
(\ref{can-lc.defn}). Finally, it can happen that several
$S_i$ come together at a point and each of them has a
  log-canonical singularity there. An explicit list is given in
\cite{ksb}. 

({\it Global property}) Instead of  each $c_1(S_i)$ being negative,
we assume that  each $c_1(S_i)-D_i$ is negative where $D_i$ is the sum of the 
double curves that lie on $S_i$. 

Another interesting issue that arises is that not every
deformation of such singular surfaces is allowed.
It turns out that even basic numerical invariants of a surface
can jump if we allow arbitrary deformations. To avoid this,
\cite{ksb} identifies a restricted deformation theory
(called $\qq$G-condition)
that produces the correct boundary.

This  answers our second Main Question:  
First, the ``simplest'' families of surfaces of general type are
families  $f:S_M\to M$ whose fibers are 
stable surfaces (and satisfy the $\qq$G-condition).
Second, every family  of surfaces of general type
is birational to such a ``simplest'' family, at least after a generically
finite-to-one change of the base $M$.

The projectivity of the resulting moduli spaces was
proved in \cite{k-proj}.
\end{say}

\begin{say}[{\bf Moduli of pairs: Alexeev, Kontsevich, 1994}]
Frequently we are interested in understanding all subvarieties $X$ of a given
variety $Y$. All is well if $X$ is smooth, but it is less clear
how to handle singular subvarieties. Various methods have been proposed,
going back to Cayley in 1860.

  Alexeev proposed in \cite{ale-pairs} that instead of working with
very singular subvarieties, one should look at
morphisms $X\to Y$ that mimic (\ref{mod.KSB}); see also 
\cite{MR2275608}. Independently, 
Kontsevich developed this approach for curves \cite{MR1291244}. 
The latter since became a standard tool in quantum cohomology theory.
\end{say}

\begin{say}[{\bf Quotient theorems: Keel, Koll\'ar, Mori, 1997}]
Step 3 of (\ref{main.steps.moduli}) leads to the general problem of
taking the quotient of a variety by a group. In our cases we have the extra
information that every point has a finite stabilizer.
In  the sixties
Artin and Seshadri proved several quotient theorems,
especially when all stabilizers are trivial.
The  general
results needed for the moduli theory were established  in 
\cite{k-quot, ke-mo}. This is a quite subtle subject since the
resulting quotients are so called {\it algebraic spaces,} a
concept somewhat more general than varieties (or even schemes).
Using the ideas of \cite{k-proj}  one can then show that, 
in the cases of interest to us, the
quotients are in fact projective (Fujino, Kov\'acs, M\textsuperscript{c}Kernan).
\end{say}

\begin{say}[{\bf Moduli in higher dimensions}]\label{mod.gen}
The  general theory follows the outlines of (\ref{mod.KSB})
with some key differences.

First, when \cite{ksb} was written, minimal models were known
to exist only in dimension 3. The higher dimensional theory
needs several results that were established only recently \cite{MR3032329}.

Second, it turned out to be quite difficult to understand how
the irreducible components of a reducible variety $X=\cup_i X_i$
glue together. For curves, as in (\ref{comp.dm.say})
the well-defined residue  of the 1-form $\frac{dx}{x}$
is a key ingredient. The current approach in higher dimension
relies on a new Poincar\'e--type residue theory for log canonical singularities;
see \cite[Chap.4]{kk-singbook}. 
The full theory should be written up in \cite{k-modbook}.
\end{say}

\begin{say}[{\bf Explicit examples: Alexeev, 2002--}]
While the above methods provide a complete answer in principle,
it has been very difficult to work out a full description in
concrete cases. The first such examples were Abelian varieties
\cite{ale-abvar} and plane curves (Hacking, \cite{hacking}). Recent
surface examples are in \cite{MR2956036}.
\end{say}

\begin{say}[{\bf Hyperbolicity: Kov\'acs, Viehweg, Zuo, 2000--2010}]
So far, very little is known about 
 the moduli spaces of surfaces and higher dimensional varieties
in general. The local structure of these spaces can be arbitrarily
complicated \cite{MR2227692}.
Hyperbolicity properties of the  moduli of smooth curves
were  conjectured by  Shafarevich in 1962 and later
extended to higher dimensions in
\cite{MR1713524, MR1978884, MR2393082, MR2726098}.
\end{say}

\begin{say}[{\bf Degeneration of Fano varieties: Xu, 2007--}]
We know much less about the moduli of Fano (=positively curved) varieties. 
Most of the geometric works deal with extending  families
  $g^*: X^*\to \Delta^*$ over a punctured disc across the puncture.
 Two questions turned out to be
especially interesting:
 understanding the combinatorial structure of
the central fiber $X_0$ for arbitrary limits
and
finding limits where  $X_0$ is especially simple.

A series of papers  \cite{k-ax, hog-xu, dkx} shows that
the combinatorial structure of $X_0$ is contractible;
this answers an old  conjecture of J.~Ax.
Recently, \cite{xu-li} shows that there are limits where $X_0$ itself is
a (singular) Fano variety, as conjectured by Tian. 
\end{say}

\begin{ack} I thank J.~Fickenscher, A.~Fulger, J.M.~Johnson, S.~Kov\'acs, 
T.~Murayama, Zs.~Patakfalvi and
N.~Sheridan
 for helpful suggestions. 
Partial financial support   was provided  by  the NSF under grant number 
DMS-0968337.
\end{ack}

%\bibliography{refs-main/refs}
\def\cprime{$'$} \def\cprime{$'$} \def\cprime{$'$} \def\cprime{$'$}
  \def\cprime{$'$} \def\cprime{$'$} \def\dbar{\leavevmode\hbox to
  0pt{\hskip.2ex \accent"16\hss}d} \def\cprime{$'$} \def\cprime{$'$}
  \def\polhk#1{\setbox0=\hbox{#1}{\ooalign{\hidewidth
  \lower1.5ex\hbox{`}\hidewidth\crcr\unhbox0}}} \def\cprime{$'$}
  \def\cprime{$'$} \def\cprime{$'$} \def\cprime{$'$}
  \def\polhk#1{\setbox0=\hbox{#1}{\ooalign{\hidewidth
  \lower1.5ex\hbox{`}\hidewidth\crcr\unhbox0}}} \def\cdprime{$''$}
  \def\cprime{$'$} \def\cprime{$'$} \def\cprime{$'$} \def\cprime{$'$}
\providecommand{\bysame}{\leavevmode\hbox to3em{\hrulefill}\thinspace}
\providecommand{\MR}{\relax\ifhmode\unskip\space\fi MR }
\providecommand{\MRhref}[2]{%
  \href{http://www.ams.org/mathscinet-getitem?mr=#1}{#2}
}
\providecommand{\href}[2]{#2}

\end{document}